# On composite conformal mapping of an annulus to a plane with two holes

Milan Batista

(July 2017)

## Abstract

In the article we consider the composite conformal map which maps annulus to infinite region with symmetric hole and nearly circular hole. It is shown that such transformation is good if the distance between centers of holes are large or radius of circular hole is small. Examples for bilinear-hypotrochoids mapping and bilinear-Schwarz-Christoffel mapping are present.

*Keywords*: Conformal mapping, composite mapping, bilinear mapping, Schwarz-Christoffel mapping

## 1 Introduction

In their recent article Lu and coworkers [1] consider a composite mapping function which maps an annulus to an infinite region with elliptic and nearly circular hole. In this short article, we will use their idea to map annulus to an infinite region with polygonal hole and nearly circular hole. To best of our knowledge such transformation is not discussed in the literature [2-9] trough it may be useful in some applications in potential theory, linear elasticity and mesh generation.

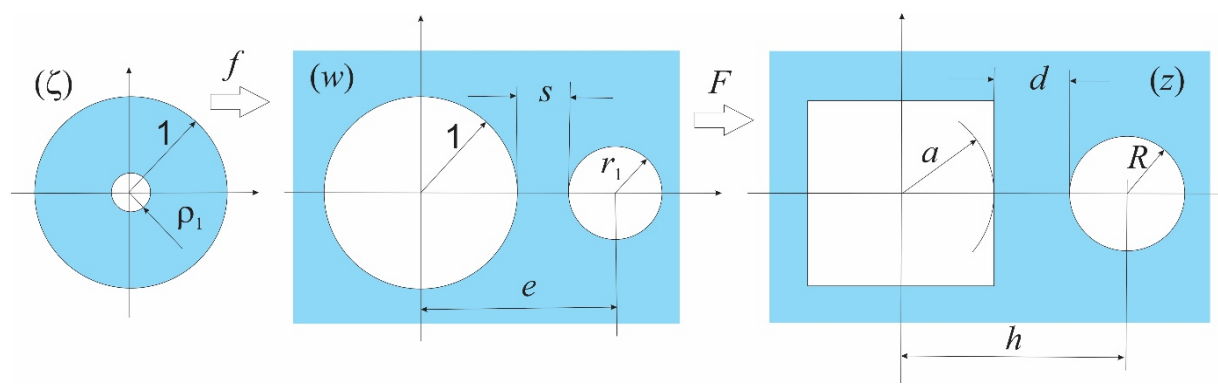

**Figure 1.** Composite mapping

## 2 The mapping

We first consider a mapping from complex plane $\zeta$ to the complex plane $w$ by the bilinear function [4, 6, 10] in the form

$$w = f(\zeta) = \frac{\zeta - \lambda}{\lambda \zeta - 1} \tag{1}$$

where $\lambda \in \mathbb{R}$. If $\lambda > 1$ and $\rho_1 < \lambda^{-1}$ then this function conformally maps concentric annulus $\rho_1 \leq |\zeta| \leq 1$ on $\zeta$ plane onto infinite regions with two non-overlapping circular holes on $w$ plane. Circle $|\zeta| = 1$ is mapped onto the circle $|w| = 1$ while circle $|w| = \rho_1$ is map onto the circle $|w - e| = r_1$ where $e_1 \in \mathbb{R}$ for non-overlapping circles obey the following inequality

$$e > 1 + r_1 \tag{2}$$

We note that $w = f(\zeta)$ has pole on annulus i.e. $f(\lambda^{-1}) = \infty$.

If $r_1$ and $e$ are given then $\rho_1$ and $\lambda$ are calculated by [10]:

$$\rho_1 = \frac{r_1 + e - \lambda}{\lambda(r_1 + e) - 1} \tag{3}$$

$$\lambda = \frac{1}{2e}\left[ e^2 + 1 - r_1^2 + \sqrt{\left(e^2 - \left(1 + r_1\right)^2\right)\left(e^2 - \left(1 - r_1\right)^2\right)} \right] \tag{4}$$

Let us now consider a function

$$z = F(w) = C\left(w + \sum_{n=1}^{\infty} \frac{c_n}{w^n}\right) \tag{5}$$

$z, w \in \mathbb{C}$, $C, c_n \in \mathbb{C}$, that is analytic in the domain $|w| > 1$. $C$ is scale factor and can be uniquely determinate from normalization e.g. $F(1) = a$, where $a$ is some characteristic length. Combining (1) and (5) we obtain

$$z = F[f(\zeta)] = C\left[ \frac{\zeta - \lambda}{\lambda \zeta - 1} + \sum_{n=1}^{\infty} c_n \left( \frac{\lambda \zeta - 1}{\zeta - \lambda} \right)^n \right] \tag{6}$$

Because $f$ and $F$ are analytic with non-zero derivatives on their domain of definition the mapping of annulus from $\zeta$ plane onto $z$ plane is conformal.

If on $z$-plane we have some simple nonintersecting curve $L$ enclosing coordinate origin then one can find coefficients of the function $F(w)$ that will map circle $|w| = 1$ to $L$ [10]. Obviously the circle $|\zeta| = 1$ on $\zeta$–plane will be by the composite mapping $F[f(\zeta)]$ mapped onto curve $L$. The circle $|\zeta| = \rho_1$ will be mapped onto the circle $|w - e| = r_1$ and $F(w)$ will map this circle onto a curve given by

$$z(\vartheta) = C\left( e + r_1 e^{i\vartheta} + \sum_{n=1}^{\infty} \frac{c_n}{\left(e + r_1 e^{i\vartheta}\right)^n} \right) \tag{7}$$

Expanding this into power series of $e^{i\vartheta}$ we obtain

$$z(\vartheta) = h + R\, e^{i\vartheta} + O\left(\varepsilon^2\right) \tag{8}$$

where we denote

$$h = C\left( e + \sum_{n=1}^{\infty} \frac{c_n}{e^n} \right) \tag{9}$$

$$R = C\left( 1 - \sum_{n=1}^{\infty} \frac{n c_n}{e^{n+1}} \right) r \tag{10}$$

$$\varepsilon = \frac{r_1}{e} = \frac{r_1}{1 + r_1 + s} \in (0,1) \tag{11}$$

$s$ is minimal spacing between the circles. Thus if $\varepsilon$ is small i.e. $r_1$ is small or $s$ is large, then the circle $|w - e| = r_1$ is approximately map onto the circle $|z - h| = |R|$.

In general, we cannot select location of circle $h$ and compute $e$, becouse $e$ is real number and therefore if $h$ is given then (9) leads to two equations for one unknown. However if curve $L$ is symmetric with respect to $x$ axis then the coefficients $c_n$ of expansion (5) are real numbers and therefore are $h$ and $R$. If in this case $h$ is given then (9) become the equation for calculating $e$

$$e + \sum_{n=1}^{\infty} \frac{c_n}{e^n} - \frac{h}{C} = 0 \qquad (12)$$

Practically we deal with finite number of terms $N$ in series. Therefore above reduce to polynomial equation of order $N+1$. However, because by assumption coefficients $c_n$ forms convergent series and because $e > 1$ the equation can be more effectively solved by some iteration method. For example, one may use direct iteration scheme in the form

$$e_0 = \frac{h}{C}$$
$$e_k = \frac{h}{C} - \sum_{n=1}^{N} \frac{c_n}{e_{k-1}^n} \qquad (13)$$

Once we know $e$ we can for given $R$ calculate $r_1$ using (10)

$$r_1 = \frac{R}{C \left( 1 - \sum_{n=1}^{N} \frac{n c_n}{e^{n+1}} \right)} \qquad (14)$$

Knowing $e$ and $r_1$ we can then calculate $\rho_1$ and $\lambda$ of bilinear mapping (1) by means of (3) and (4). In this way we obtain a mapping which conformal maps annulus to infinite region bounded by the curve $L$ and nearly circular hole $|z - h| = R$.

## 3 Examples

### 3.1 Composite bilinear-hypotrochoids mapping

Consider the mapping [10]

$$z = F(w) = C \left( w + \frac{m}{w^n} \right) \qquad (15)$$

where $0 \le |m| \le \frac{1}{n}$, and $n \in \mathbb{N}$. By this function the circle $|w| = 1$ is mapped onto a curve which is bounded in a region

$$r_{in} = C(1 - m) \le \left| F\left( e^{i\vartheta} \right) \right| \le C(1 + m) = r_{out} \qquad (16)$$

The curve minimal radius of curvature at $\vartheta_{\min} = 0$ and maximal radius of curvature at $\vartheta_{\max} = \pi/(n+1)$ are

$$r_{\min} = C\frac{(mn-1)^2}{mn^2+1} \tag{17}$$

$$r_{\max} = C\frac{(mn+1)^2}{mn^2-1} \tag{18}$$

If $r_{\max} = \infty$ then

$$m = 1/n^2 \tag{19}$$

In words, the figure has locally straight edges. In Figures 2, 3, 4, 5 are shown some of these maps.

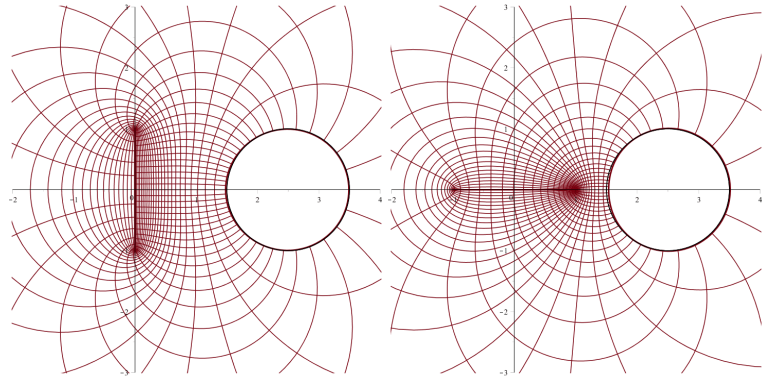

**Figure 2.** Composite bilinear-hypotrochoids mapping for the case $n=1$, $d=0.5$, $r_{out}=1$ and $r=1$. The maximal discrepancies from circle are 0.0211 and 0.0353.

Let us now consider a discrepancy of the mapping to the given circle. Using (9) and (10) for $e$ and $r_1$ can write the difference between the circle $|z-h| = R$ and the mapping curve as

$$\Delta = f\left(e + r_1 \mathrm{e}^{i\vartheta}\right) - \left(h + R\mathrm{e}^{i\vartheta}\right) = \frac{C}{n^2 e^n}\left[n\varepsilon\mathrm{e}^{i\vartheta} - 1 + \frac{1}{\left(1+\varepsilon\mathrm{e}^{i\vartheta}\right)^n}\right] \tag{20}$$

For $\varepsilon \ll 1$ we have

$$\Delta = \frac{C}{n^2 e^n} \frac{\varepsilon^2}{2}\left(1+\frac{1}{n}\right)e^{2i\vartheta} + O\left(\varepsilon^3\right) \tag{21}$$

As was already shown the difference becomes smaller with decreasing $\varepsilon$, however it also decrease by increasing $n$.

Consider now a case when $s \to 0$ i.e. the circles on $w$ plane almost touch. In this case by using (11) for $\varepsilon$ we can write (20) as

$$\Delta = \frac{C}{n^2}\left[\frac{\left(ne^{i\vartheta}-1\right)r_1-1}{\left(1+r_1\right)^n} + \frac{1}{\left(1+r_1\left(1+e^{i\vartheta}\right)\right)^n}\right] + O\left(s_1\right) \tag{22}$$

At $e^{i\vartheta} = -1$ we obtain maximum difference

$$\Delta_{max} = \frac{C}{n^2}\left[1 - \frac{\left(n+1\right)r_1+1}{\left(1+r_1\right)^n}\right] + O\left(s_1\right) \tag{23}$$

So

$$\lim_{r_1 \to \infty} \Delta_{max} = \frac{C}{n^2} \tag{24}$$

i.e. the difference is bounded and decrease with increasing $n$. This property can be also observed from numbers in Table 1.

**Table 1.** Difference $\Delta_{max}$ for the case $m = 1/n^2$, $n = 2$, $r_{out} = 1$, and different size of hole $R$ and spacing $d$. ($\varepsilon = \frac{r_1}{e}$)

| $d$ | $10^{-5}$ | | $0.1$ | | $1$ | |
|---|---|---|---|---|---|---|
| $R$ | $\varepsilon$ | $\Delta_{max}$ | $\varepsilon$ | $\Delta_{max}$ | $\varepsilon$ | $\Delta_{max}$ |
| 0.25 | 0.2600 | 0.0294 | 0.2248 | 0.0170 | 0.1151 | 0.0012 |
| 0.5 | 0.3804 | 0.0522 | 0.3474 | 0.0350 | 0.2051 | 0.0036 |
| 1 | 0.5261 | 0.0794 | 0.4974 | 0.0587 | 0.3382 | 0.0087 |
| 2 | 0.6764 | 0.1033 | 0.6537 | 0.0810 | 0.5030 | 0.0164 |
| 4 | 0.8025 | 0.1181 | 0.7866 | 0.0956 | 0.6679 | 0.0240 |
| 8 | 0.8894 | 0.1247 | 0.8796 | 0.1023 | 0.8003 | 0.0288 |
| 128 | 0.9922 | 0.1280 | 0.9915 | 0.1058 | 0.9846 | 0.0320 |

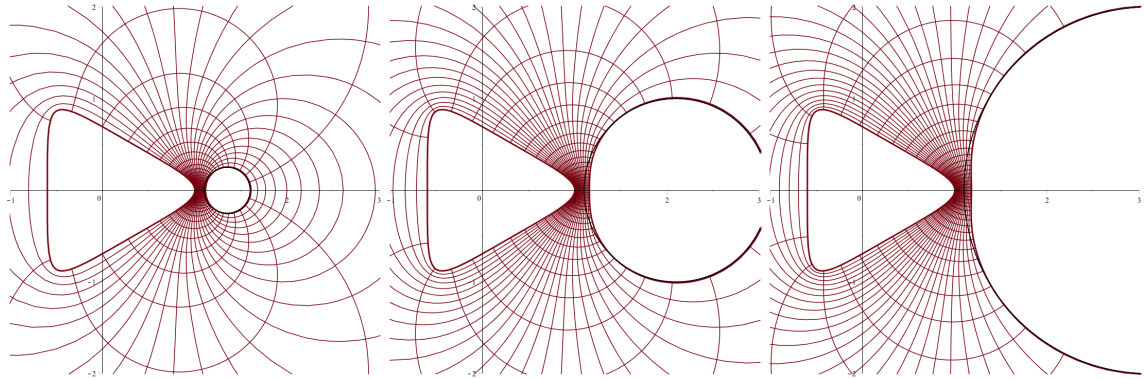

**Figure 3.** Composite bilinear-hypotrochoids mapping. Effect of hole size for $m = 1/n^2$, $n = 2$, $r_{out} = 1$, $d = 0.1$ $r \in \left\{ \frac{1}{4}, 1, 2 \right\}$. The maximal differences from the circle are 0.0170, 0.0587, 0.0810.

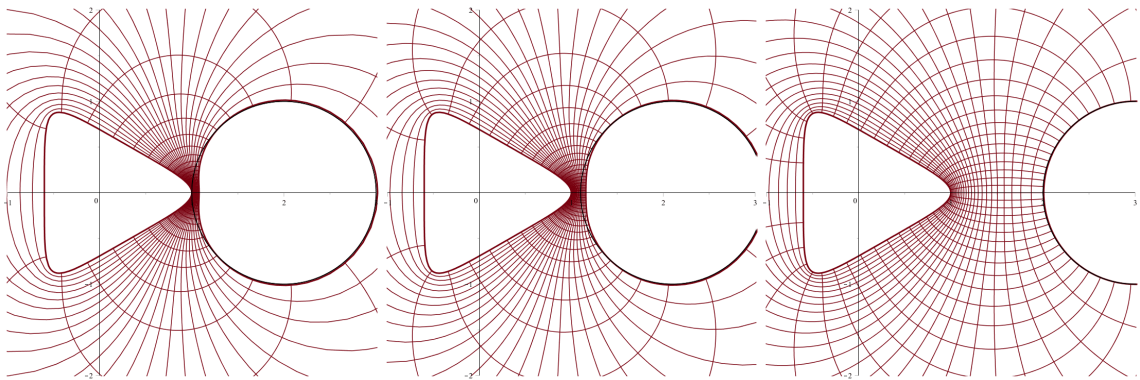

**Figure 4.** Composite bilinear-hypotrochoids mapping. Effect of hole position for $m = 1/n^2$, $n = 2$, $r_{out} = 1$, $r = 1$ and $d \in \left\{ 10^{-5}, 0.1, 1 \right\}$. The maximal differences from the circle are 0.0794, 0.0587 and 0.0087

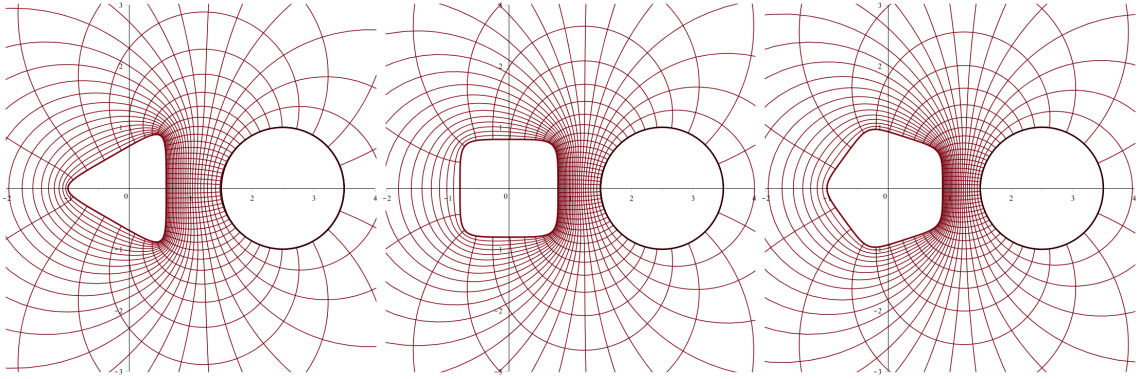

**Figure 5.** Composite bilinear-hypotrochoids mapping. Effect of hole shape for $m = -1/n^2$, $d = 1/2$, $r_{out} = 1$, $r = 1$ and $n \in \{2, 3, 4\}$. The maximal differences from the circle are 0.0184, 0.0110, 0.0060.

### 3.2 Bilinear - Schwarz–Christoffel mapping

For regular polygon of $n$ sides the Schwarz-Christoffel maping has the form [11]

$$z = F(w) = C \int \left(1 - \frac{1}{w^n}\right)^{\frac{2}{n}} d\zeta \tag{25}$$

The integral can be evaluated by series expansion

$$z = F(w) = C \sum_{k=0}^{\infty} (-1)^k \binom{2/n}{k} \frac{w^{-nk+1}}{1-nk} \tag{26}$$

For a figure rotated by angle $\frac{2\pi}{n}$ the mapping is slightly different

$$z = F(w) = C \sum_{k=0}^{\infty} \binom{2/n}{k} \frac{w^{-nk+1}}{1-nk} \tag{27}$$

In Figures 6 and 7 are show some of these maps.

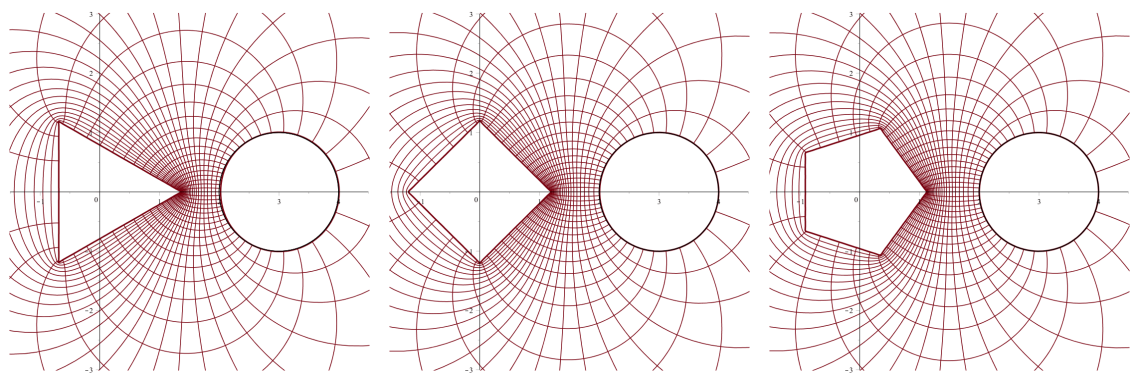

**Figure 6.** Composite bilinear-Schwarz-Christoffel mapping with five terms. Effect of shape for $C = 1$, $d = 1$ and and $r = 1$ for $n \in \{3, 4, 5\}$. The maximal differences from the circle are 0.02524, 0.0089, 0.0034.

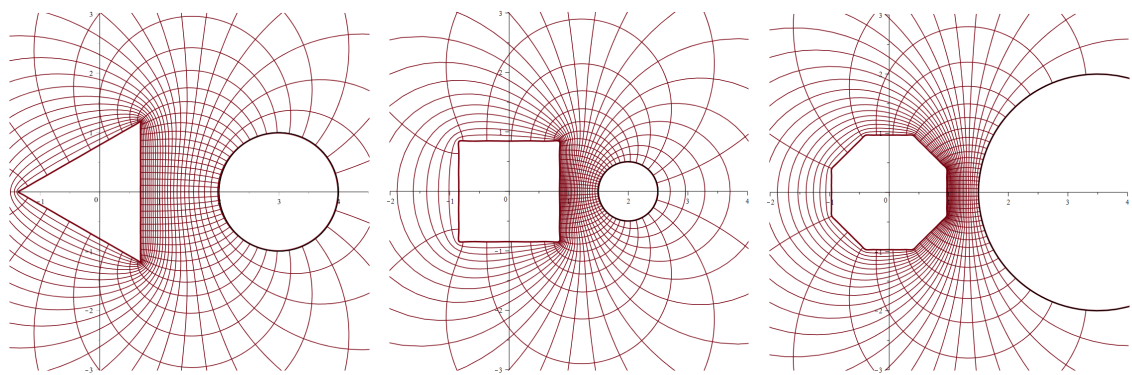

**Figure 7.** Composite bilinear-Schwarz-Christoffel mapping with three terms. For triangular hole $C = 1$, $d = 1$, $R = 1$, for square hole $a = 1$, $d = 0.5$, $R = 0.5$, and for octagonal hole $a = 1$, $d = 0.5$, $R = 2$. The maximal differences from the circle are 0.0188, 0.0182, 0.0271